\documentclass[11pt]{article}
\usepackage{amssymb}
\usepackage{amsmath}
\usepackage{amsfonts}
\usepackage[a4paper]{geometry}

\newtheorem{defi}{Def\hskip 1pt inition }[section]
\newtheorem{lem}[defi]{Lemma}
\newtheorem{theo}[defi]{Theorem}
\newtheorem{prop}[defi]{Proposition}

\newtheorem{assum}[defi]{Assumption}
\newcommand{\comment}[1]{}

\def\videbox{\mathbin{\vbox{\hrule\hbox{\vrule height1ex \kern.5em\vrule
height1ex}\hrule}}}
\def\1{1\!{\rm l}}

\def\tlam{{\Lambda_*}}

\def\A{{\mathcal A}}
\def\G{{\mathcal G}}
\def\D{{\mathcal D}}
\def\E{{\mathbb E}}
\def\Var{{\rm Var}}
\def\L{{\mathbb L}}
\def\N{{\mathbb N}}
\def\P{{\mathbb P}}
\def\R{{\mathbb R}}

\begin{document}

\title{Approximation of additive random fields based on standard information:
average case and probabilistic settings}

\author{Mikhail Lifshits
\footnote{Dept. Mathematics and Mechanics, St.Petersburg State University, 198504, Bibliotechny pr. 2,  Russia, and MAI, 
Link\"oping University, Sweden}
\and 
Marguerite Zani
\footnote{Laboratoire d'Analyse et Math\'ematiques appliqu\'ees , UMR CNRS 8050, Universit\'e Paris-Est -Cr\'eteil, 61 av du G\'en\'eral de Gaulle, 94010, Cr\'eteil Cedex, France
} 
} 

\maketitle

\abstract{
We consider approximation problems for tensor product and additive random fields based on standard 
information in the average case setting.  We also study the probabilistic setting of the mentioned
problem for tensor products.
The main question we are concerned with in this paper is {``How much do we loose by considering standard 
information algorithms against those using general linear information?''}
For both types of the fields, the error of linear algorithms has been studied in great detail. However, the power of
standard information for them was not addressed so far, which we do here. Our main conclusion is that in most 
interesting cases there is no more than a logarithmic loss in approximation error when information is 
being restricted to the standard one.

The results are obtained by randomization techniques.
}

\section{Introduction: general information against standard information}

Let $f$ be a function which we consider as an element of some Banach space $(B,||\cdot||)$ of 
functions. Assume that the whole function $f$ is unknown but we are able to measure the 
values of some functionals $F_1(f),...,F_n(f)$, such as values at certain points, 
integrals, etc. Let $\Psi:\R^n\to B$ be some mapping. Then we may call
\[
  Af:= \Psi(F_1(f),...,F_n(f))
\] 
an {\it approximation algorithm} and evaluate its error by $||f-Af||$. 
A typical problem setting assumes minimization of approximation error
for given $n$ by optimizing the choice of the set $(F_j(\cdot))_{1\le j\le n}$
within available class of functionals, as well as the choice of the mapping $\Psi$.

As for the target function $f$, we essentially have two options for problem setting. 
Either we can let it vary over some set, typically some compact
subset of $B$, for example, the unit ball of some other Banach space compactly embedded
into $B$, by taking the {\it worst possible} approximation error for algorithm evaluation.
Or we may consider $f$ as a {\it random} function having certain distribution in $B$ and using
the expectation of the error for algorithm evaluation. The two possibilities are often referred 
to as the "worst case setting", resp. "average case setting". We will rather stick to the
latter one and consider random functions (or {\it random fields}, if their arguments are multivariate) 
as the objects of approximation. 

Let us now specify the problem furthermore. We will use a Hilbert norm for error evaluation,
which essentially means  $B=\L_2([0,1]^d)$, where $d$ is understood as appropriate parametric dimension
of a random field we consider. It is well known that for Hilbert norms one may restrict the
class of all algoritms to the class of {\it linear algorithms}, which means
\[
   A f = \sum_{j=1}^n F_j(f) \phi_j, 
\]
where every $F_j(\cdot)$ is a linear functional on $B$ and $\phi_j\in B$ are specific fixed elements. 
It is often said that such algorithms are based on {\it linear information}. 
In many applications it is too costly to calculate values of arbitrary linear functionals needed
for such approximation. It is therefore preferable to restrict the choice of $F_j(\cdot)$ by letting 
$F_j(f):=f(t_j)$, i.e. by taking values of $f$ at some points. One says that such algorithms are based  
on {\it standard information}.

It is clear that functionals  based on linear information form a larger class than those based on standard 
information, thus, after optimization, the average error would be smaller in linear case. 
The main question we are concerned with in this paper is
{\it ``How much do we loose by considering standard information algorithms against those using general linear 
information?''}

This problem already received much interest, see \cite{HWW,HNV,KWW,WW2} and especially the survey \cite{NW_survey}.
Most of the research was concentrated on the worst case setting. It was shown that in many cases the polynomial
term of error decay (when considered as a function of the number $n$ of functionals used) is the same for
algorithms based on linear and standard information. However, there exist cases where the behavior of the error as 
a function of $n$ for linear and standard information is drastically different. It means that the problem
we consider is by far non-trivial. 

In this paper we restrict investigation to one specific but very important class of averaging function distributions,
or, equivalently, random fields under consideration. Namely we consider {\it tensor product random fields} and 
{\it additive random fields}. In short,
tensor product random field  on $[0,1]^d$ is a zero mean random field with covariance function
\[
    K_d(s,t)  :=  \prod_{l=1}^d K(s_l,t_l)  
\]
where $s=(s_1,\cdots,s_d), t=(t_1,\cdots,t_d)\in [0,1]^d$,
and $K(\cdot,\cdot)$ is a covariance function on  $[0,1]^2$. 
Brownian sheet, Brownian pillow, and other famous random fields 
belong to this class.

Next, a $d$-parameter additive random field of order $b$, defined on
$[0,1]^d$, appears as follows.
We pick any set of $b$ coordinates from $d$ available ones, consider
a tensor product random field introduced above but depending only of the picked coordinates
and then sum up uncorrelated versions of those tensor products over all possible sets of
coordinates of size $b$. As a simplest example, if $b=1$, we just have a random field
\[
   X(t) = \sum_{l=1}^{d} X_l(t_l), 
\]
where $X_l(s), 0\le s\le 1$, are uncorrelated copies of a one-parametric field with covariance
$K(\cdot,\cdot)$. See more definition details in Sections \ref{s:tensor} and  \ref{s:add}, respectively.

For both types of the fields, the error of linear algorithms is studied in great detail,
cf. \cite{LT,LPW, LZ, WW}, as well as closely related small deviation behavior \cite{KNN}. However, the power of
standard information for them was not addressed so far, which we do here. Our main conclusion is that in most 
interesting cases there is no more than a logarithmic loss in approximation error when information is 
being restricted to the standard one.

Notice that we do not merely consider a fixed random field but also explore the case 
when one or even two parameters, dimension $d$ and additivity order $b$ go to infinity, in which case
a concept of {\it relative} error is relevant. 
 
There is still much room for the research in the direction we traced here. We do not really pretend
that approximation rates for standard information case obtained in this work are optimal, although 
they are obviosly quite close to those ones. Neither we provide concrete algorithms related to obtained rates. 
In fact, our results are rather existence theorems because they are obtained by randomization techniques.
These interesting topics still wait for subsequent, much deeper investigation.
\bigskip

The paper is organized as follows. In Section \ref{s:algo} we recall a quite 
general approximation procedure providing upper bounds for the quality of algorithms 
based on standard information assuming that similar bounds are known for algorithms 
based on general linear functionals. This procedure is systematically applied in the 
sequel to different settings and to different random fields of interest. 

In Section \ref{s:tensor} we recall the notion of a tensor product random field and 
the known results about the quality of their approximation based on general linear 
functionals. Furthermore, we obtain new results for approximation of those fields 
based on standard information. Consideration concerns both fixed and increasing 
parametric dimension of random field.

In Section \ref{s:add} we perform the same program for more complicated additive 
random fields and come out with the new estimates  for approximation of those fields 
based on standard information.  

Finally, in Section \ref{s:prob} we come back to tensor product random fields in 
fixed dimention and evaluate the power of standard information in the probabilistic 
setting, i.e. searching for the algorithms that assure given approximation quality 
with prescribed probability.

\section{Approximation procedure based on standard information}
\label{s:algo}

\subsection{The pointwise approximation algorithm}

Let $(D,\nu)$ be a measure space. We first consider an approximation procedure for a deterministic
function $g\in \L_2(D,\nu)$ based on the standard information, i.e. on the values of $g$ at $n$ 
points $\tau_1,\dots, \tau_n$.

The randomized algorithm below is a slight refinement of those presented in the series of works  \cite{HWW}, \cite{HW},
 \cite{KWW}, and \cite{WW2}.

Fix two integers $n,m$. Let $(\eta_j)_{j\geq 1}$ be an orthonormal system in  $\L_2(D,\nu)$. Then
\[
    u_m(t):=\frac{1}{m}\sum_{j=1}^m\eta_j^2(t)
\]
is a density on $D$. Given $n$ points $\tau=(\tau_1,\cdots,\tau_n)\in D^n$ we define a linear approximation algorithm
$A_{\tau}: \L_2(D,\nu)\to \L_2(D,\nu)$ as follows. Let $g\in \L_2(D,\nu)$. Then $g$ admits the Fourier
expansion
$$
  g=\sum_{j=1}^m g_j\eta_j+g_m^{\perp}, 
$$
where 
$$g_j=\langle g,\eta_j\rangle:=\int_D g\eta_j d\nu \,,$$
and $g_m^{\perp}$ is orthogonal to the space spanned by $\{\eta_j,j\leq m\}$.
Let 
\[
    A_{\tau} g :=\sum_{j=1}^m \hat g_j\eta_j\,,
\]
where
$$
    \hat g_j := \frac{1}{n}\sum_{l=1}^n g (\tau_l) \, \frac{\eta_j(\tau_l)}{u(\tau_l)} 
$$
is a pointwise approximation to the integral $g_j$.

Consider $\tau$ as an i.i.d. sample in $D$ with density $u_m$. For any $j,l$ we can compute
\[
       \E_{\tau}\left( g(\tau_l)\, \frac{\eta_j(\tau_l)}{u_m(\tau_l)}\right)
       =\int g(t)\, \frac{\eta_j(t)}{u_m(t)} \, u_m(t) \, \nu(dt) = g_j \,.
\]
Hence, $\E_{\tau}(\hat g_j)=g_j$. 

We can evaluate the variance as well,

\[
   \mbox{\Var}_{\tau}\left(g(\tau_l) \, \frac{\eta_j(\tau_l)}{u_m(\tau_l)}\right)
   \leq 
   \E_{\tau} \left[\left(g(\tau_l) \, \frac{\eta_j(\tau_l)} {u_m(\tau_l)}\right)^2\right]
   =
   \int g(t)^2  \, \frac{\eta_j(t)^2}{u_m(t)}  \, \nu(dt)  \,.
\]
Hence,
\[
     \mbox{\Var}_{\tau}(\hat g_j)
     =
     \E_{\tau}[(\hat g_j-g_j)^2]
     \leq 
     \frac{1}{n} \int g(t)^2 \, \frac{\eta_j(t)^2}{u_m(t)}  \, \nu(dt)   \,.
\]
By summing up,

\[
   \E_{\tau}\left[\sum_{j=1}^m(\hat g_j-g_j)^2\right]
   \leq 
   \frac{1}{n} \int g(t)^2  \, \frac{\sum_{j=1}^m   \eta_j(t)^2}{u_m(t)}\, \nu(dt)
   =
   \frac{m}{n} \int g^2(t)  \, \nu(dt)
   = \frac{m}{n} \, \|g\|^2_2
\,.
\]

Recall that 
\[
    \|g-A_{\tau}g\|_2^2=\sum_{j=1}^m(\hat g_j-g_j)^2+\|g_m^{\perp} \|_2^2   \, ,
\]
therefore,

\begin{equation}\label{pyth0}
     \E_{\tau}\|g-A_{\tau}g\|_2^2
     \leq
     \frac{m}{n}  \, \|g\|_2^2+\|g_m^{\perp}\|_2^2   \, .
\end{equation}

Now we pass to approximation of {\it random} functions.
Let us consider a random field $Y(t,\omega)$ on $D$ and, assuming the sample
$\tau$ independent of $Y$, apply inequality  
(\ref{pyth0}) to the sample paths $g(t):=Y(t,\omega)$. 
By Fubini theorem,
\[
   \E_{\tau} \E_Y  \|Y-A_{\tau}Y\|_2^2
   =
   \E_Y  \E_{\tau}  \|Y-A_{\tau}Y\|_2^2
   \leq
   \frac{m}{n}  \, \E \|Y\|_2^2 + \E \|Y^{\perp}_m\|_2^2,
\]
where $Y^{\perp}_m$ is the orthogonal part of $Y$ with respect to the finite 
system  $\{\eta_j\}_{j=1}^{m}$,

\[
   Y(t,\omega)=\sum_{j=1}^m  (Y,\eta_j)  \, \eta_j(t)+Y^{\perp}_m(t,\omega)   \,.
\]

By the mean value theorem, there exists a point set  
\[
   \tau:=\tau(Y):=\{\tau_l(Y),l=1\cdots m\}
\]
such that 

\begin{equation}\label{pyth2}
    \E_{Y}  \|Y-A_{\tau}Y\|_2^2   
    \leq  
    \frac{m}{n} \, \E\|Y\|_2^2+ \E \|Y_m^{\perp}\|_2^2   \,.
\end{equation}
We stress the fact that the set $\tau(Y)$ does not depend on the sample paths of $Y$,  
i.e. it is non-random.

\subsection{Iteration procedure}
\label{ss:iterA}

Now we define an iteration procedure giving a sequence of approximations $A_k, k\geq 0,$ for a field $Y$. 
Notice that the parameters $n,m$ and the orthonormal system $(\eta_j)_{j=1}^m$ are fixed throughout 
the construction. 
Let $A_0:=0$ and iteratively
\[
      A_k Y :=  A_{k-1} Y +A_{\tau(Y-A_{k-1}Y)} (Y-A_{k-1} Y ) \, ,
\]
Let denote $\Delta_k:=Y-A_k Y$.  Then by construction,
\begin{eqnarray*}
    \Delta_k &=&  Y-A_{k-1} Y - A_{\tau(Y-A_{k-1}Y)}  (Y-A_{k-1}Y)
\\
             &=&\Delta_{k-1}-A_{\tau(\Delta_{k-1})}(\Delta_{k-1}) \,  .
\end{eqnarray*}

By applying inequality (\ref{pyth2}) to the field $\Delta_{k-1}$ instead of $Y$, we get

\[
     \E \|\Delta_k\|_2^2    
     \leq     
     \frac{m}{n}  \,   \E\|\Delta_{k-1}\|_2^2+  \E  \|\Delta_{k-1}^{\perp}\|_2^2 \,  .
\]

Notice that any algorithm used keeps invariant the orthogonal part, hence
 $$
    \E \|\Delta_{k-1}^{\perp} \|_2^2 =    \E \|Y_m^{\perp}\|_2^2\,.
 $$
 Then we obtain an iterative estimate 
 \begin{equation}\label{err}
        \mbox{err}_k    \leq \frac{m}{n}     \,    \mbox{err}_{k-1}+R,
 \end{equation}
where $\mbox{err}_k:=\E\|\Delta_k \|_2^2$ and $R:=\E\|Y_m^{\perp}\|_2^2$ does not depend on $k$.

The iterative sequence like this quickly approaches a limit assuming that
$
\tfrac{m}{n}<1 \,  .
$

Indeed, let $x:=R/(1-\frac{m}{n})$ and write (\ref{err}) as

\[
     \mbox{err}_k   \leq   \frac{m}{n}  \  \mbox{err}_{k-1}+(1-\frac{m}{n}) \, x
\]
which is equivalent to
\[
     \mbox{err}_k-x\leq\left(\frac{m}{n}\right)^k(\mbox{err}_0-x),
\]
whence
\[
     \mbox{err}_k\leq x+\left(\frac{m}{n}\right)^k\mbox{err}_0.
\]
By the previous definitions, the latter inequality reads as
\begin{equation}\label{pyth3}
   \E \|\Delta_k\|_2^2  
   \leq   
   \frac{\E\|Y_m^{\perp }\|_2^2}  
         {\scriptstyle 1-\frac{m}{n}}+
   \left(\frac{m}{n}\right)^k \E \|Y\|_2^2  \, .
\end{equation}
Notice that the algorithm $A_k$ uses $k n$ point evaluations ($n$ new points on each step).

Now we adjust the algorithm by connecting the variables $m$ and $n$. 

For any $n$ take $m:=[n/2]$ and apply algorithm  $A_{n,k}:=A_k$. We get by (\ref{pyth3})
\begin{equation}\label{pyth4}
  \E \|Y-A_{n,k} Y \|_2^2  
  \leq 
  2 \E\|Y_{[n/2]}^{\perp}\|^2_2 +  2^{-k}E\|Y\|^2_2 \,  .
\end{equation}

Let us now make an assumption about the average approximation rate of $Y$ by finite parts 
of its Fourier expansion. We assume the following {\it polynomial approximation rate}:

\begin{assum} \label{a:polynomial}
There exist $p>0, C_1>0$, and a real $\gamma$ such that
for any $m\geq 1$,
\begin{equation}  \label{rate_fourier}
   \E \Big\|Y- \sum_{j=1}^m (Y,\eta_j) \eta_j \Big\|^2_2   
   = \E \|Y_m^{\perp}\|^2_2  \leq C_1  m^{-2p}    \,  (\ln m)^\gamma \,  \E\|Y\|^2_2\,  . 
\end{equation} 
\end{assum}

By letting $k:= [Z\log_2 n]+1$ with large $Z>2p$, we get the folowing

\begin{eqnarray} \nonumber 
   \E\|Y-A_{n,k} Y\|_2^2
   &\leq& 
   2 C_1 [n/2]^{-2p} (\ln [n/2])^\gamma   \E\|Y\|^2_2 +   n^{-Z} \E\|Y\|^2_2
   \\  \label{polerr} 
   &\leq& 
   C n^{-2p}   (\ln n)^\gamma \,  \E\|Y\|^2_2,
\end{eqnarray}
while $([Z \log_2 n]+1)n$ point evaluations of $Y$ are used and $C=C(C_1,p,b,Z)$.
\medskip

We may summarize those findings in the following result.

\begin{prop} 
Let $Y(t), t\in D$, be a random field satisfying Assumption \ref{a:polynomial} with some parameters $p,\gamma, C_1$.
Let $Z>2p$. Then there exists $C=C(C_1,p,\gamma,Z)$ such that for every positive integer $n$ there exist
a linear approximation procedure $A_n$ based on $([Z \log_2 n]+1) \times n$ point evaluations of $Y$ 
and satisfying
\begin{equation} \label{rate_point}
   \E\|Y-A_{n} Y\|_2^2 
   \leq  
   C n^{-2p}   (\ln n)^\gamma \,  \E\|Y\|^2_2. 
\end{equation}

\end{prop}

{\bf Remark.}
Notice that the approximation rate while using point measurements \eqref{rate_point} 
is the same, up to a constant, as in the case of using arbitrary functionals
 \eqref{rate_fourier}. The minor loss consists in logarithmic increase of the number of used measurements:
 $([Z \log_2 n]+1) n$ instead of $n$. Unfortunately, our result  does not provide construction of the measurement 
points because of the use of randomization technique.

\section{Tensor product random fields}
\label{s:tensor}

\subsection{Simple tensor product}
\label{ss:tensor}

Let us apply this previous scheme to {\it tensor product} random fields. First, recall the latter notion.
Let $K(\cdot,\cdot)$ be a covariance kernel on $[0,1]$. By the well known Mercer's theorem it can 
be represented as
\[
    K(s,t)  =  \sum_{i=1}^\infty \lambda(i)^2 \phi_{i}(s) \phi_{i}(t)
\]
where $(\phi_i)_{i>0}$ is an orthonormal system of functions in $\L_2([0,1])$, and  $(\lambda(i)^2)_{i>0}$
is the sequence of eigenvalues (corresponding to the integral operator with kernel $K$) satisfying
condition
\[
    \sum_{i=1}^\infty \lambda(i)^2 <\infty.
\]

We consider a random field $X(t), t\in [0,1]^d$, with zero mean and covariance function
\begin{equation} \label{Ktensor}
    K_d(s,t)  :=  \prod_{l=1}^d K(s_l,t_l)  
\end{equation}
where $s=(s_1,\cdots,s_d), t=(t_1,\cdots,t_d)\in [0,1]^d$. It is natural to call $X$ the $d$-th
tensor degree of a one-parametric random field with covariance $K$.

It is well known that  $X$ admits Karhunen--Lo\`eve expansion

\begin{equation} \label{tensor}
    X(t)   =   \sum_{k\in\N^d}\prod_{l=1}^d\lambda(k_l)\prod_{l=1}^d\phi_{k_l}(t_l) \, \xi_k \
    :=  \sum_{k\in\N^d} \lambda_k \phi_k \, \xi_k, 
\end{equation}
where 
$\{\xi_k,k\in\N^d\}$ are non-correlated random variables with zero mean and unit variance;
and
\[
  \lambda_k :=\prod_{l=1}^d\lambda(k_l), 
  \quad
  \phi_k(t) :=\prod_{l=1}^d\phi_{k_l}(t_l)   \  ,
  \quad
  k\in\N^d   \,.
\]
The sample paths of $X$ are considered as elements of the space $\L_2(D,\nu)$ where $D=[0,1]^d$ and
$\nu$ is the $d$-dimensional Lebesgue measure.

In the paper \cite{LT}, Lifshits and Tulyakova studied approximation rate for 
such random fields, considering fixed and increasing dimension $d$ for the case when the measurment of arbitrary
linear functionals is available. In the present work, we aim to explore the power of approximation  
algorithms based upon {\it pointwise} evaluations of $X$. In terms of information based complexity theory \cite{NW,TWW, W},
we consider standard information.

\subsubsection{Fixed dimension}
\label{ssec_tensor_fd}
For fixed $d$ we assume

\begin{assum}\label{eigen0}
\[
   \lambda(i)\sim \mu \, i^{-r}(\ln i)^q\,, i\to\infty
\]
for some $\mu>0$, $r>1/2$ and $q\neq r$.
\end{assum}

There are many important examples of such behavior, including Wiener process, fractional Brownian motion, etc, see
e.g. \cite{L12}.

Let $(\bar\lambda_j^2,j\in\N)$ be the decreasing rearrangement of the array $(\lambda_k^2\,,k\in\N^d)$. 
In fixed dimension, for a sequence of eigenvalues satisfying Assumption \ref{eigen0},  we have the following elementary result
for which we refer to \cite{LT}:

\begin{lem}\label{lem}
Let $\alpha :=q/r$. Then
\begin{eqnarray}
\label{eigen1}
   \overline\lambda_{j}^2 
   &\sim & B_d^2\,  j^{-2r}(\ln j)^{2r\beta}\,,\quad j\to\infty,
   \\
   \nonumber
   \\
   \nonumber
   \mbox{where }
   &\bullet& \textrm{for } \alpha>-1:
   \begin{cases}
   B_d=\mu^d \left(\frac{\Gamma(\alpha+1)^d}{\Gamma(d(\alpha+1))}\right)^r , \\
   \beta=(d-1)+d\alpha , 
   \end{cases}
\\
    \nonumber
    &\bullet& \textrm{for } \alpha<-1: 
    \begin{cases}
       B_d= \mu d^r \left[\sum_{i\geq 1}\lambda(i)^{1/r}\right]^{(d-1)r}  , \\
       \beta=\alpha.
    \end{cases}
\end{eqnarray}
\end{lem}
Equivalent results can be found e.g. 
in Cs\'aki \cite{Cs1}, Li  \cite{Li1},  Papageorgiou and Wasilkowski \cite{PW} (for $q=0$) 
and in Karol', Nikitin, and Nazarov \cite{KNN}.

The asymptotics for $\alpha=-1$ is also known; it has a form
\[
 \overline\lambda_{j}^2 
   \sim  B_d^2\,  j^{-2r}(\ln j)^{-2r}(\ln\ln j)^{2r(d-1)} \,,\quad j\to\infty,
\]
cf. \cite{KNN}. We skip this case here and in the sequel just for the sake of brevity.

By summing up and using \eqref{rate_point}, we have the following polynomial approximation error:
\begin{equation}
\label{eigensum}
      \E \|X-X_m\|_2^2 =   \sum_{j>m}   \overline   \lambda_j^2 
      \leq  C_1 m^{1-2r}  (\ln m)^{2r\beta}
\end{equation}
with some $C_1=C_1(X)$.

Let us consider the expansion (\ref{tensor}) and apply the preceding procedure
with the orthonormal system $(\bar\phi_j,j\in\N)$ being the 
rearrangement of the system  $(\phi_k)$  in (\ref{tensor}) corresponding to the 
decreasing order of eigenvalues $\lambda_k^2$.

We obtain from (\ref{polerr}) that for every $n$ there is a pointwise algorithm $A_n$
that uses  $([Z \log_2 n] +1) n$ point evaluations and gives an error

\[
   \E \|X-A_n X \|_2^2
   \leq C n^{1-2r}   (\ln n)^{2r \beta} 
\]
with some  $C=C(X)$ not depending of $n$.

A trivial variable change leads to the final result:

\begin{prop} \label{prop23}
Let $X$ be a tensor product random field $(\ref{tensor})$.
Under Assumption $\ref{eigen0}$ there exists a constant $C=C(X)$ and a sequence 
of linear approximation algorithms 
$(A_n)_{n\ge 1}$ using respectively $n$ point values such that  for each $n$ we have
\begin{equation}  \label{polerr3} 
       \E \| X-A_{n} X \|_2^2  \leq C \ n^{1-2r}   (\ln n)^{2r (\beta+1)-1}. 
\end{equation}
\end{prop}

\subsubsection{Increasing dimension}
\label{increasord1}

Now we let $d$ go to infinity and consider $d$-parametric random fields as usually done in
information based complexity theory \cite{NW, TWW, W}. We still consider tensor product random fields  \eqref{tensor}
with covariance \eqref{Ktensor} generated by a fixed univariate covariance function $K(\cdot,\cdot)$ but
denote the field $X_d$ instead of $X$, in order to stress the dependence of parametric dimension.

In increasing dimension, one must consider {\it relative} error, since the total 
range of the random field
\begin{equation}  \label{Lambdad}
    \E\|X_d\|_2^2=\left( \sum_{i=1}^{\infty} \lambda(i)^2\right)^d
    :=  \Lambda^d
\end{equation}  
may go to infinity as $\Lambda>1$,  resp. to zero as $\Lambda<1$, whenever the dimension $d$ increases.
Therefore, for a given $\epsilon\in (0,1)$ we look for a pointwise approximation algorithm $A$ such that 
\[
   \E \|X_d-A X_d\|_2^2 \leq \epsilon^2\Lambda^d\, ,
\]
and ask how many point measurements of $X_d$ we need for achieving this.

Let the cardinality associated to the 
{\it relative} error  $\tilde m_d(\epsilon)$ be defined by

\begin{equation*}
     \tilde m_d(\epsilon)=  \inf   \left\{   
     m \,  ; \    \sum_{j>m}   \overline   \lambda_j^2  
     \leq   \epsilon^2\Lambda^d  
     \right \}\, ,
\end{equation*}
From \cite{LT} we know that the following result holds.

\begin{theo} $(\cite{LT})$ \label{thm_LT2} 
Under assumption 
\begin{equation}\label{M2}   
   M_2:=   \sum_{i=1}^{\infty}|\ln \lambda(i)|^2\lambda(i)^2<\infty\,,
\end{equation}
for any $\epsilon\in (0,1)$ it is true that
\[
     \lim_{d\to\infty} \frac{\ln\tilde m_d(\epsilon)-d\ln \tlam}{\sqrt d}=2q_*\,,
\]
where 
\[
   M:=-\sum_{i=1}^{\infty}\ln\lambda(i)\frac{\lambda(i)^2}{\Lambda},
\] 
\begin{equation}  \label{tlam} 
    \tlam:=\Lambda e^{2M},
\end{equation}
and the quantile $q_*$ is given by the equation
$
   \hat \Phi\left(\frac{q_*}{\sigma}\right)=\epsilon^2\, ,
$
where $\sigma^2=\frac{M_2}{\Lambda}-M^2$ and
$$
\hat\Phi(x)=\frac{1}{\sqrt{2\pi}}\int_x^\infty \exp\{-u^2/2\}du
$$
is the tail of the standard normal law.
\end{theo}

With pointwise estimation, from \eqref{pyth4} and \eqref{Lambdad} we have with 
$n=2 \tilde m_d(\epsilon)$

\[
   \E \|   X_d-A_{n,k}  X_d \|_2^2\leq (2\epsilon^2+2^{-k})\Lambda^d \, .
\]

Taking $k=[-\log_2(1/\epsilon)^2)]+1$ points, we have the relative error $3\epsilon^2$
while using only  $2 ([-\log_2(\epsilon^2)]+1)  \tilde m_d(\epsilon) $ points.

Our conclusion is that the number of necessary measurements explodes almost in the same way,  
whenever we use general or pointwise functionals. Using notations from Theorem
 \ref{thm_LT2} we have again

\begin{prop}   
Under assumption \eqref{M2}
for any $\epsilon\in (0,1)$ there exists a sequence of positive integers $m_d(\epsilon)$
and a sequence of linear approximation procedures $A_d^\epsilon$ based upon
$m_d(\epsilon)$ point measurements of $X_d$ such that for any $d$
\[
   \E \|   X_d-A_{d}^\epsilon  X_d \|_2^2 \leq 3\epsilon^2  \E \|X_d\|_2^2 \, .
\]
and  
\[
     \lim_{d\to\infty} \frac{\ln\tilde m_d(\epsilon)-d\ln \tlam}{\sqrt d}=2q_*\,,
\]
 
\end{prop}

\section{Additive random fields}
\label{s:add}

\subsection{Additive random field of order $b$}

We recall now the definition of a $d$-parameter additive random field of order $b$, cf.
\cite{CL}, \cite{LZ}.
Essentially it means  that we pick any set of $b$ coordinates from $d$ available ones, consider
a tensor product random field introduced in Subsection \ref{ss:tensor} 
above but depending only of the picked coordinates
and then sum up uncorrelated versions of those tensor products over all possible sets of
coordinates of size $b$. 

By the reasons that will soon become clear, this time it is convinient for us to numerate
the eigenfunctions and eigenvalues of the tensor product starting from zero. Thus let denote
$\N_0:=\{0,1,2,\dots\}$ in addition to $\N:=\{1,2,\dots\}$. 

Let $\lambda(i)\ge 0,  i\in\N_0$, satisfy $\sum_{i=0}^\infty \lambda(i)^2<\infty$  
and let $(\phi_i)_{i\in\N_0}$ be an orthonormal basis in $\L_2([0,1])$.

Let denote by $\D$ and $\D_b$ the following sets of indexes:
\[
   \D=\{1,\cdots, d\}\,,\quad \D_b=\{\A\subset \D\,,|\A|=b\}\,.
\]
Here and elsewhere $|\A|$ stands for the number of elements in a finite set $\A$.
For each $\A\in \D_b$ we write $\A=\{a_1,\cdots ,a_b\}$ with elements ordered, say,
as $a_1<\dots< a_d$.

Now an additive $d$-parametric random field $X_{d,b}(t)$, $t\in[0,1]^d$,
of order $b$ is defined by
\begin{equation}
 \label{tensor_add}
       X_{d,b}(t) :=  \sum_{\A\in \D_b} \sum_{k\in\N_0^\A}\left(\prod_{a\in \A}\lambda(k_a) 
       \prod_{a\in \A}\phi_{k_a}(t_a)\right)\xi_{k}^\A\,.
\end{equation}
Here $(\xi_{k}^\A,\ \A\in \D_b, k\in \N_0^\A)$ is a family of non-correlated random 
variables with mean zero and unit variance.

We assume furthermore
\begin{assum}\label{zero}
$$
    \forall u\in[0,1]\,,\,\,\phi_0(u)=1\,.
$$
\end{assum}
This technical assumption ensures (see \cite{LZ}, Lemma 2.2) that the family 
\[
    \Big\{ \phi_{k}^\A(t):= \prod_{a\in \A}\phi_{k_a}(t_a), \ \A\subset \D, k\in\N^\A \Big\} 
\]
is an orthonormal system in $\L_2([0,1]^d)$. It appears in \eqref{tensor_add}.
Many interesting and important processes satisfy Assumption \ref{zero}, see \cite{HWW}.

Under this assumption,  \eqref{tensor_add} transforms into 
\begin{equation}
 \label{tensor_add2}
       X_{d,b}(t) =  \sum_{\A\in \D_b}  \sum_{\G\subset\A} \lambda(0)^{|\A|-|\G|} \sum_{k\in\N^\G}
       \left(\prod_{a\in \G}\lambda(k_a) \right) \phi_{k}^\G(t) \xi_{k_0}^\A\, ,
\end{equation}
where $k_0\in (\N_0)^\A$ is obtained from $k\in \N^\G$ by adding zeros. 

Looking at the expansion \eqref{tensor_add2}, we see that the same eigenfunction
$\phi_{k}^\G(t)$ appears as many times as there are sets $\A\in\D_b$ containing $\G$.
Since the number of such sets is $C_{d-|\G|}^{b-|\G|}$, we may write 
\begin{eqnarray}
       \nonumber
       X_{d,b}(t) &=&  \sum_{{\G\subset \D\atop |\G|\le b}}  \left(C_{d-|\G|}^{b-|\G|}\right)^{1/2} \lambda(0)^{b-|\G|} 
       \sum_{k\in\N^\G}
       \left(\prod_{a\in \G}\lambda(k_a) \right) \phi_{k}^\G(t) \widetilde\xi_{k}^\G 
       \\ \label{tensor_add3}
        &=& \sum_{h=0}^b  \sum_{\G\in \D_h}  \left(C_{d-h}^{b-h}\right)^{1/2} \lambda(0)^{b-h} 
       \sum_{k\in\N^\G}
       \left(\prod_{a\in \G}\lambda(k_a) \right) \phi_{k}^\G(t) \widetilde\xi_{k}^\G \, ,
\end{eqnarray}
where $(\widetilde\xi_{k}^\G,\ \G\subset \D, |\G|\le b, k\in \N^\G)$, is a family of non-correlated random 
variables with mean zero and unit variance.

\subsubsection{Fixed dimensions}

In this subsection, we assume that dimension $d$ and order $b \le d$ are fixed, and for an integer $h$ satisfying 
$0\le h\le b$ denote by
\[
    \lambda_{k,h}^2=\prod_{l=1}^h\lambda(k_l)^2\,,\qquad k=(k_1,\cdots, k_h)\in \N^h.
\]    
Looking at \eqref{tensor_add3}, we see that eigenvalues of covariance operator have the form
\[
  C_{d-h}^{b-h}\lambda(0)^{2(b-h)} \lambda_{k,h}^2 
\]
with multiplicity $|\D_h|=C_d^h$.

For every fixed $h$ let denote by $(\overline \lambda_{j,h}^2\,,j\in\N)$ the decreasing rearrangement of the 
array $(\lambda_{k,h}^2, k\in \N^h)$.

As in case of fixed dimension for simple tensor product (Subsection \ref{ssec_tensor_fd}), 
let Assumption \ref{eigen0} on the asymptotics of $\lambda(i)$ hold with some parameters
$r$ and $q$. 
Then by Lemma \ref{lem} the asymptotics of 
$\overline\lambda_{j,h}^2$ is given by \eqref{eigen1}.

Subsequent analysis leads to different conclusions in two cases.

\begin{itemize}
\item
If $\alpha:=q/r >-1$, then the exponent $\beta=h-1+h\alpha$  in  (\ref{eigen1}) 
depends on $h$. We observe the slowest decay of $\overline\lambda_{j,h}^2$
for $h$ delivering maximal $\beta$, i.e. for $h=b$. 

When we merge all arrays and rearrange the total set of eigenvalues, the asymptotics will be the same as
for the one for the dominating array.  It follows that for appropriately chosen orthonormal system 
$\overline\phi_j$ in $\L_2([0,1]^d)$ we have the Fourier expansion
\begin{equation}\label{decompo2}
      X_{d,b}=\sum_{j=1}^{\infty} (X_{d,b}, \overline\phi_j) \, \overline\phi_j 
\end{equation}
with the error 
\[
  X_{m}^{\perp} := X- \sum_{j=1}^{m} (X_{d,b}, \overline\phi_j)\, \overline\phi_j 
  =\sum_{j=m+1}^{\infty} (X_{d,b}, \overline\phi_j) \, \overline\phi_j 
\]
admitting the bound
\[
   \E\|X_{m}^{\perp}\|_2^2 \le C m^{1-2r}  (\ln m)^{2r\beta} 
\]
as in \eqref{eigensum} with $\beta=b-1+b\alpha$.

Thus we arrive at the same conclusion as in the case of simple tensor field, i.e. Proposition \ref{prop23}. 
Namely, we have

\begin{prop} 
Let $X_{d,b}$ be an additive tensor product random field of order $b$ defined in $(\ref{tensor_add})$.
Under Assumption $\ref{eigen0}$ (with $q>-r$) and Assumption $\ref{zero}$ there exists a constant 
$C=C(X_{d,b}$) and a sequence of linear approximation algorithms 
$(A_n)_{n\ge 1}$ using respectively $n$ point values such that for each $n$ we have
\[
       \E \|   X_{d,b}-A_{n}  X_{d,b} \|_2^2  \leq C\ n^{1-2r}   (\ln n)^{2b(r+q)-1}. 
\]
\end{prop}

The proof is exactly the same as for Proposition \ref{prop23}. Just notice the calculation
of logarithm's exponent: recall that $h=b$, $\alpha=q/r$,  and $\beta=h-1+h\alpha$
yield  for the exponent in \eqref{polerr3} 
\[
    2r(\beta+1)-1 =2r h(1+\alpha) -1= 2b(r+q)-1.
\]  

\item
If $\alpha=q/r<-1$, then $\beta=\alpha$ does not depend on $h$ when Lemma \ref{lem} is applied. 
The rearranged eigenvalues $\overline\lambda_{m,h}^2$ have the same order of decay for all $h$,
namely,
\[
  \overline\lambda_{m,h}^2 \, \sim B_h^2 \, j^{-2r}\, (\ln j)^{2q}
\]
where the constants $B_h$ are defined in Lemma \ref{lem} above.
Moreover, the same decay order corresponds to the rearranged total set of eigenvalues.

Recall how the asymptotically optimal finite rank approximation based on arbitrary functionals was constructed in
\cite{LZ}. Let us construct an optimal $m$-term approximation. Let denote by 
$m_h$ the number of largest eigenvalues in the array $\overline \lambda_{j,h}^2$ 
contributing to the approximation. We know from \cite{LZ} that the quasi-optimal choice is
\[
   m_h:= \left[ m\ \frac{Q(h)^{1/2r}}{\sum_{l=1}^b  Q(l)^{1/2r}} \right] \,,
\]
where, taking into account the eigenvalue factor 
$C_{d-h}^{b-h}\lambda(0)^{2(b-h)}$ that was left out when moving from the true eigenvalues
to  $\lambda_{k,h}^2$, and the multiplicity $C_d^h$, 
\[
     Q(h):= C_{d-h}^{b-h}\, \lambda(0)^{2(b-h)}  \, B_h^2  \, [C_d^h]^{2r}  ,
\]
it is easy to see that the error is of order
\begin{eqnarray*}
   &&\sum_{h=1}^b   C_{d-h}^{b-h} \, \lambda(0)^{2(b-h)} \  C_d^h  \sum_{j>m_h/C_d^h}  B_h^2 \, j^{-2r} \, (\ln j)^{2q}
   \\
   &\sim&
      (2r-1)^{-1}  \sum_{h=1}^b   C_{d-h}^{b-h} \, \lambda(0)^{2(b-h)} \  C_d^h  B_h^2 \,  
      \left(m_h/C_d^h\right)^{1-2r}  \, (\ln m_h)^{2q}
   \\
   &=&
    (2r-1)^{-1}   \sum_{h=1}^b   Q(h)  m_h^{1-2r} \, (\ln m_h)^{2q}
   \\ 
   &\sim&
    (2r-1)^{-1}   \sum_{h=1}^b   Q(h) \ Q(h)^{(1-2r)/2r} \left(\sum_{l=1}^b  Q(l)^{1/2r}\right)^{2r-1}  
     m^{1-2r} \, (\ln m)^{2q}
   \\
    &=& Q\,  m^{1-2r}(\ln m)^{2q},
\end{eqnarray*}
where
\[
    Q := (2r-1)^{-1} \left( \sum_{h=1}^b  Q(h)^{\frac{1}{2r}} \right)^{2r}.
\]
This means that an ordering of orthogonal basis is possible 
as in expansion (\ref{decompo2}) with the error  of $m$-term approximation
\[ 
  \E\|X_{m}^{\perp}\|_2^2 \sim  Q\,  m^{1-2r}(\ln m)^{2q}.
\]
Now we have the following.

\begin{prop} 
Let $X_{d,b}$ be an additive tensor product random field of order $b$ defined in $(\ref{tensor_add})$.
Under Assumption $\ref{eigen0}$ (with $q<-r$) and Assumption $\ref{zero}$ there exists a constant 
$C=C(X_{d,b}$) and a sequence of linear approximation algorithms 
$(A_n)_{n\ge 1}$ using respectively $n$ point values such that for each $n$ we have
\[
       \E \|   X_{d,b}-A_{n}  X_{d,b} \|_2^2  \leq C\ n^{1-2r}   (\ln n)^{2(q+r)-1}. 
\]
\end{prop}

The proof is again the same as for Proposition \ref{prop23}, while the calculation
of logarithm's exponent is as follows.

The relations $\beta=\alpha=q/r$ yield  for the exponent in \eqref{polerr3} 
\[
    2r(\beta+1)-1 =2r\beta +2r-1= 2q+2r-1.
\]

\end{itemize}

\subsubsection{Increasing dimension}

As in Subsection \ref{increasord1},
we let parameter dimension $d$ go to infinity and consider $d$-parametric random fields as usually done in
information based complexity theory \cite{NW}. In this context of growing $d$ we 
consider  additive random fields $X_{d,b}$ of order $b$ described
in  \eqref{tensor_add}. Recall that in increasing dimension, one must consider {\it relative} error, since the total 
range of the random field, measured by its mean square varies with dimension as
\[
   \E(\|X_{d,b}\|^2_{\L_2([0,1]^d}) = C_d^b \, \Lambda^b\, ,
\]
where  $ \Lambda:= \sum_{i=0}^{\infty} \lambda(i)^2$.

Recall that the spectrum of the process is described as follows. To any fixed $h=1,\dots,b$
associate an array of eigenvalues
\begin{equation}
    \label{eig_arr}
    \Big\{ 
         \lambda_{k,h}^2 :=\prod_{l=1}^h\lambda(k_l)^2\,,\qquad k=(k_1,\cdots, k_h)\in \N^h.
     \Big\}
\end{equation}
Then multiply all eigenvalues by  $C_{d-h}^{b-h}\, \lambda(0)^{2(b-h)}$ and take them with
multiplicity $C_d^h$.

Therefore, for given $\epsilon$,
we look for a minimal $n$ such that an algorithm $A_n$ based on $n$ point evaluations of $X_{d,b}$
provides 
\[
      \E\|X_{d,b}-A_n\, X_{d,b}\|_{\L_2([0,1]^d}^2 \leq C_d^b \, \Lambda^d \, \epsilon^2 .
\]

Admitting that dimension $d$ goes to infinity, we have a choice between keeping the additivity order
$b$ fixed or letting it go to infinity, too. Consider both cases.

\paragraph{Order $b$ fixed and dimension $d$ increasing.}
Let  $n_{d,b}^{avr}(\epsilon)$ be the minimal $m$ such that a representation 
\begin{equation}\label{decompo3}
      X_{d,b}=\sum_{j=1}^{m} ( X_{d,b}, \phi_j) \, \phi_j+X^{\perp}\, ,
\end{equation}
holds with some orthogonal system $(\phi_j)_{1\le j\le m}$
and
\[
   \E||X^{\perp}||^2_2\le C_d^b  \Lambda^b \epsilon^2. 
\]
It is known from \cite{LZ} that if Assummption \ref{eigen0} and Assumption
\ref{zero} hold, if $\epsilon\to 0$, $d\to\infty$, and 
$b$ is fixed,  that the main contribution 
to the approximation is given by the eigenvalue array \eqref{eig_arr} corresponding to $h=b$. Therefore,
\[
   n_{d,b}^{avr}(\epsilon)  \sim \frac{d^b}{b!}\ \Lambda^{-b/(2r-1)} \ n_b^{avr}(\epsilon)
\]
where the approximation cardinality $n_b^{avr}(\epsilon)$ is defined in \cite{LT}, 
for a simple tensor product random field in dimension $b$ (additive field of order $b=d$) and
behaves as follows
\[
   n_b^{avr}(\epsilon) 
   \sim
   \left(\frac{B_b}{\sqrt 2(r-1/2)^{r\beta+1/2}}
   \frac{|\log\epsilon|^{r\beta}}{\epsilon}\right)^{1/(r-1/2)}\,, \qquad \epsilon\to 0,
\]
with $B_b$ and $\beta$ given in Lemma \ref{lem} (let $d=b$ in those notations) and depending on
$b$ and on parameters $r,q$ from Assumption \ref{eigen0} but not depending on $d$. 
Notice especially that $ n_{d,b}^{avr}(\epsilon)$ depends on $d$ polynomially.

Now we transform the approximation algorithm \eqref{decompo3} that uses general linear functionals
$(\phi_j,\cdot)$ into a point evaluation based algorithm $A_{n,k}$ from Section \ref{ss:iterA}.
We still use parameters $n=2 n_{d,b}^{avr}(\epsilon)$ and $k=[\log_2 (1/\epsilon)]+1$ and obtain 
from \eqref{pyth4}
\begin{equation}\label{pyth4bd}
  \E \|X_{d,b}-A_{n,k} X_{d,b} \|_2^2  
  \leq 
  3 C_d^b  \Lambda^b \epsilon^2  ,
\end{equation}
whyle using 
\[
  n\ k = 2\  n_{d,b}^{avr}(\epsilon) \, ([\log_2(1/\epsilon)]+1)
\]
point measurements.
Thus the polynomial dependence of number of points on $d$ remains the same. On the other hand, 
the number of points keeps its polynomial term in dependence of $\epsilon$ while getting an extra degree
of logarithmic term. The relative error slightly grows up from $\epsilon$  to $\sqrt{3}\epsilon$.   
One can, of course, arrange a trade off between the constants in the growth 
of relative error and number of used points.

\paragraph{Order $b$ and dimension $d$ increasing proportionally.}

In this case we assume that both $b$ and $d$ tend to infinity while $b/d\to f\in(0,1)$. 
It was observed in \cite{LZ} that the main 
contribution to the field range $\E\|X_{d,b}\|_2^2$ is given by the eigenvalue arrays 
\eqref{eig_arr} corresponding to $h$ such that $h/b\sim p$ where 
$p=1-\tfrac{\lambda(0)^2}{\Lambda}$. 
Using this observation, it was shown that for any fixed $\epsilon\in(0,1)$  
and assuming \eqref{M2} along with Assumption \ref{zero}, 
the representation (\ref{decompo3}) is possible
with the same error $E||X^{\perp}||^2_2\le C_d^b  \Lambda^b \epsilon^2$ 
while the dimension $n$ grows with a certain exponential rate. More precisely,
for each small $\delta$ and all large $d>d_0(\delta)$ we have 
$n\le V^{(1+\delta)d},$ 
where
\[
    V:=(1- f p)^{f p-1} f^{-f p}(1-p)^{(1-p)f}\tlam^{f}\, .
\]
and $\tlam$ depending on eigenvalues $(\lambda(i))$ was defined in \eqref{tlam}

Recall some special cases.
\begin{itemize}
 \item If $f=1$, then 
\[ 
     V=(1-p)^{p-1} (1-p)^{1-p}\tlam =\tlam. 
\]
This essentially corresponds to the case considered in Section \ref{increasord1}.
  \item It is rather surprising that that for $f<1$ the explosion coefficient $V$
depends on the ``strange'' quotient $p$. 
\item If $p=1$, i.e. $\lambda(0)=0$, then we have a nice formula
\[
    V=(1-f)^{f-1} f^{-f}\tlam^{f}\, .
\]
\item If $f=0$, then $V=1$, which means that there is no exponential explosion. This includes the case
of $b$ fixed while $d$ growing to infinity.  
\end{itemize}

By applying the pointwise algorithm $A_{n,k}$ with the same $k=[\log_2(1/\epsilon)]+1$ iterations we get
again the bound \eqref{pyth4bd}.  
The number of points used 
\[
  n\ k \le 2 V^{(1+\delta)d} ([\log_2(1/\epsilon)]+1)
\]
is still exponential in variable $d$. On the exponential scale, the explosion coefficient $V$ 
remains the same when we pass to pointwise algorithms.

\section{Probabilistic setting}
\label{s:prob}



In this section, we consider a simple tensor product random field $X(t,\omega)$ 
\begin{equation} \label{tensor2}
    X(t,\omega)   =   \sum_{k\in\N^d}\prod_{l=1}^d\lambda(k_l)\prod_{l=1}^d\phi_{k_l}(t_l)\xi_k(\omega) \
   :=  \sum_{k\in\N^d} \lambda_k \phi_k \xi_k,
\end{equation}
(exactly as in formula \eqref{tensor}, with the notations explained there).
But since we will deal now with probabilistic estimates, we assume additionaly that random variables
$\xi_k$ are {\it jointly Gaussian}. Since $\xi_k$ are supposed to be non-correlated, in Gaussian context
this implies their independence.
We also suppose that Assumption \ref{eigen0} is verified.

For given $\epsilon>0$ and $\gamma\in (0,1)$, we are interested in a linear algorithm $A$ based on point 
measurements of $X$, satifying
\begin{equation} \label{epsgamma}
  \P(  \|X-A X\|_2>\epsilon) \leq\gamma\,
\end{equation}
and using the minimal possible number of point values.

\begin{prop} \label{p:prob}
Let $X$ be a tensor product random field $\eqref{tensor2}$
with eigenvalues satisfying Assumption $\ref{eigen0}$.
Let  $\epsilon>0$ and $\gamma\in (0,1)$. Denote $v:=\epsilon\left( 1+\sqrt{2|\ln \gamma|} \right)^{-1} $.
Then for some constant $C=C(X)$ there exist 
a linear approximation algorithm 
$A$ such that \eqref{epsgamma} holds 
and using at most
\begin{equation} \label{nv}
    C \left(v^2 |\ln v|^{1-2r(\beta+1)} \right)^{-1/(2r-1)}
\end{equation}
point values, where $r$ comes from Assumption $\ref{eigen0}$ and $\beta$ is given in Lemma $\ref{lem}$.
\end{prop}


For proving Proposition \ref{p:prob}, we will need the following simple fact from the 
theory of Gaussian vectors.

\begin{lem} \label{l:gauss}
Let $X$ be a zero mean Gaussian vector in a separable Banach space. Then for any $\gamma\in (0,1)$
we have 
\[
   \P\left( ||X||\ge \left[ \E||X||^2\right]^{1/2} \left( 1+\sqrt{2|\ln \gamma|} \right) \right)
   \le \gamma.
\]
 \end{lem}

\noindent {\bf Proof of Lemma \ref{l:gauss}.} 
Let us write  $||X||$ in the dual form,
\[
   ||X||=\sup_{||f||=1} (f,X),
\]
where supremum is taken over the unit sphere of the dual space.
Notice that $(f,X)$ is a centered Gaussian random variable and
\[
   \E(f,X)^2 \le \E (||f||\, ||X||)^2 =  ||f||^2\, \E ||X||^2.  
\]
It follows that
\[
  \sigma^2 := \sup_{||f||=1}  \E(f,X)^2 \le  \E ||X||^2.
\] 
We will also use the known fact that a median $med(||X||)$ satisfies
inequality
\[
   med(||X||) \le \E ||X|| \le \left(\E ||X||^2\right)^{1/2},
\]
see \cite{L}. Now by Gaussian isometric inequality, see \cite{L} again,
\begin{eqnarray*}
   \P\left( ||X||\ge \left[ \E||X||^2\right]^{1/2} \left( 1+\sqrt{2|\ln \gamma|} \right) \right)
   &\le&
  \P\left( ||X|| \ge  med(||X||)  + \sigma \sqrt{2 |\ln \gamma|} \right) 
  \\
    &\le& 1- \Phi\left(\sqrt{2|\ln \gamma|} \right)
   \\
    &\le& \exp\left( - \left(\sqrt{2|\ln \gamma|} \right)^2/2\right)=\gamma.
\end{eqnarray*}
Here $\Phi(\cdot)$ is the distribution function of the standard normal law.
$\Box$
\medskip

\noindent {\bf Proof of Proposition \ref{p:prob}.} According to \eqref{polerr3}, for every $n$
we may find an algorithm $A_n$ using $n$ point values such that
\[   
       \E \| X-A_{n} X \|_2^2  \leq C_0 \ n^{1-2r}   (\ln n)^{2r (\beta+1)-1}. 
\]
Now choose a minimal $n$ so large that
\[  
     C_0 \  n^{1-2r}   (\ln n)^{2r (\beta+1)-1} \le v^2
\]
which amounts to say that $n$  exceeds expression \eqref{nv} 
with appropriate constant $C(X)$ independent of $n$. Let $A:=A_n$.
Then 
\[
  \left[ \E||X-A X \|_2^2 \right]^{1/2} \left( 1+\sqrt{2|\ln \gamma|} \right)
  \le v \left( 1+\sqrt{2|\ln \gamma|} \right) = \epsilon
\]
and by applying Lemma \ref{l:gauss} to the Gaussian vector $X-A X$ we get
\[
   \P(  \|X-A X\|_2>\epsilon) \leq
 \P\left( ||X-A X||\ge \left[ \E||X-A X||^2\right]^{1/2} \left( 1+\sqrt{2|\ln \gamma|} \right) \right)
\le
\gamma.
\]
$\Box$

\bibliographystyle{plain}

\small

\begin{thebibliography}{10}

{\baselineskip=12pt

\bibitem{CL}
 Chen, X., Li, W.V.(2003) Small deviation estimates for some additive processes, In: 
``High dimensional probability, III (Sandjberg, 2002)'', Ser.: Progr. Probab., 55, 
Birkh\"auser, Basel, 225--238. 

\bibitem{Cs1}
Cs\'aki, E. (1982) On small values of the square integral of a multiparameter
Wiener process, In: "Statistics and Probability. Proc. III Pannonian Symp.
 Math. Statist.", D.\,Reidel, Boston, 19--26.

\bibitem{HW}
Hickernell F.J., Wo\'zniakowski, H. (2000)  Integration and approximation in  arbitrary
dimension, Adv. Comput. Math., 12, 25--58.

\bibitem{HWW} 
Hickernell, F.J., Wasilkowski, G.W., and Wo\'zniakowski, H. (2006) 
Tractability of linear multivariate problems in the average-case setting.
Proc. Conf. Monte-Carlo and Quasi-Monte Carlo Methods. Ulm, 2007, 
(Heinrich, S., Keller, A., Niederreiter, H., eds.), Springer, 461--493.

\bibitem{HNV} Hinrichs, A.,  Novak, E., and Vybiral, J. (2008)  
Linear information versus function evaluations for L2-approximation,
J. Approx. Theory, 153, 97--107. 

\bibitem{KNN} Karol', A.I., Nazarov, A.I., and Nikitin, Ya.Yu. (2008)
Small ball probabilities for Gaussian random fields and 
tensor products of compact operators,
Trans. Amer. Math. Soc.,  360, 1443--1474. 

\bibitem{KWW} Kuo, F., Wasilkowski, G., and Wo\'zniakowski, H.  (2009)
On the power of standard information for multivariate approximation in the worst case setting,
J. Approx. Theory, 158, 97--125.

\bibitem{Li1}
Li, W.V. (1992) Comparison results for the lower tail of Gaussian
seminorms, J. Theor. Probab., 5, 1--31.

\bibitem{L}
Lifshits, M.A. Gaussian Random Functions. Kluwer, 1996.

\bibitem{L12}
Lifshits, M.A. Lectures on Gaussian Processes. Springer, 2012.

\bibitem{LPW}
Lifshits, M.A., Papageorgiou A., and Wo\'zniakowski, H. (2012) Average case 
tractability of non-homogeneous tensor product problems, J. Complexity, 
28, 539--561.

\bibitem{LT}
Lifshits, M.A., Tulyakova, E.V. (2006) Curse of dimensionality in approximation 
of Gaussian random fields, Probab. Math. Statist., 26, no. 1, 97--112.

\bibitem{LZ}
Lifshits, M.A., Zani, M. (2008) Approximation complexity of additive random fields,
J. Complexity, {24}, no. 3, 362--379.

\bibitem{NW_survey}
Novak, E.,  Wo\'zniakowski, H. (2011)
On the power of function values for the approximation problem in various settings,
Surveys in Approx. Theory, 6, 1-23. 

\bibitem{NW}
Novak, E.,  Wo\'zniakowski, H. Tractability of Multivariate Problems,  
volumes 1--3, EMS, Z\"urich, 2008-2012. 

\bibitem{PW} Papageorgiou, A., Wasilkowski, G.W. (1990) On the average
complexity of multivariate problems, J. Complexity, 6, 1--23.

\bibitem{TWW}
Traub, J.F., Wasilkowski, G.W., and Wo\'zniakowski, H. Information-based
Complexity, Academic Press, Boston, 1988. 

\bibitem{WW} Wasilkowski, G.W., Wo\'zniakowski, H. (1999) Weighted
tensor-product algorithms for linear multivariate problems, J. Complexity,
15, 1--56.

\bibitem{WW2} Wasilkowski, G.W., Wo\'zniakowski, H. (2007) 
The power of standard information for multivariate approximation in the randomized setting, 
Math. Comp., 76, 965--988.

\bibitem{W}
Wo\'zniakowski, H. (1994) Tractability and strong tractability of linear multivariate
problems, J. Complexity, 10, 96--128.


}
\end{thebibliography}

\end{document}